\documentclass[10pt,a4paper]{amsart}

\usepackage{amsmath}
\usepackage{amssymb}
\usepackage{amsthm}
\usepackage{amsfonts}
\usepackage{microtype}
\usepackage{mathrsfs}
\usepackage{graphicx}
\usepackage{color}
\usepackage{latexsym}
\usepackage{rotating}
\usepackage{xspace}
\usepackage[all]{xy}
\usepackage{longtable}
\usepackage{leftidx}
\usepackage{mathtools}
\usepackage{titlesec}
\usepackage{tikz}
\usetikzlibrary{calc}
\usetikzlibrary{arrows}
\usetikzlibrary{decorations.markings}
\usepackage{geometry}

\newtheorem{theorem}{Theorem}[section]
\numberwithin{equation}{theorem}
\newtheorem{lemma}[theorem]{Lemma}
\newtheorem{proposition}[theorem]{Proposition}
\newtheorem{corollary}[theorem]{Corollary}

\theoremstyle{definition}
\newtheorem{definition}[theorem]{Definition}
\newtheorem{example}[theorem]{Example}

\theoremstyle{conjecture}

\newcommand{\Ass}{\operatorname{Ass}}
\newcommand{\im}{\operatorname{im}}

\newcommand{\Spec}{\operatorname{Spec}}

\newcommand{\ara}{\operatorname{ara}}

\newcommand{\ann}{\operatorname{ann}}

\newcommand{\cd}{\operatorname{cd}}

\newcommand{\hd}{\operatorname{hd}}

\newcommand{\Cone}{\operatorname{Cone}}

\newcommand{\Ext}{\operatorname{Ext}}
\newcommand{\Supp}{\operatorname{Supp}}

\newcommand{\Tor}{\operatorname{Tor}}
\newcommand{\Hom}{\operatorname{Hom}}

\newcommand{\depth}{\operatorname{depth}}
\newcommand{\width}{\operatorname{width}}
\newcommand{\Coass}{\operatorname{Coass}}
\newcommand{\coker}{\operatorname{coker}}

\newcommand{\fm}{\frak{m}}
\newcommand{\fp}{\frak{p}}

\newcommand{\suchthat}{\;\ifnum\currentgrouptype=16 \middle\fi|\;}

\newenvironment{prf}[1][Proof]{\begin{proof}[\bf #1]}{\end{proof}}

\makeatletter
\newcommand{\holim@}[2]{%
  \vtop{\m@th\ialign{##\cr
    \hfil$#1\operator@font holim$\hfil\cr
    \noalign{\nointerlineskip\kern1.5\ex@}#2\cr
    \noalign{\nointerlineskip\kern-\ex@}\cr}}%
}
\newcommand{\holim}{%
  \mathop{\mathpalette\holim@{\rightarrowfill@\textstyle}}\nmlimits@
}
\makeatother

\makeatletter
\def\@secnumfont{\bfseries}
\def\section{\@startsection{section}{1}%
  \z@{.7\linespacing\@plus\linespacing}{.5\linespacing}%
  {\normalfont\Large\bfseries\filcenter}}
\def\subsection{\@startsection{subsection}{2}%
  \z@{.5\linespacing\@plus.7\linespacing}{-.5em}%
  {\normalfont\large\bfseries}}
\makeatother

\begin{document}

\sloppy

\author[K. Divaani-Aazar, H. Faridian and M. Tousi]{Kamran Divaani-Aazar, Hossein Faridian
and Massoud Tousi}

\title[Local Homology, Koszul Homology and Serre Classes]
{\Large Local Homology, Koszul Homology and Serre Classes}

\address{K. Divaani-Aazar, Department of Mathematics, Alzahra University, Vanak, Post Code 19834, Tehran,
Iran-and-School of Mathematics, Institute for Research in Fundamental Sciences (IPM), P.O. Box 19395-5746,
Tehran, Iran.}
\email{kdivaani@ipm.ir}

\address{H. Faridian, Department of Mathematics, Shahid Beheshti University, G. C., Evin, Tehran, Iran, Zip Code 1983963113.}
\email{h.faridian@yahoo.com}

\address{M. Tousi, Department of Mathematics, Shahid Beheshti University, G.C., Tehran, Iran, P.O. Box 19395-5746,
Tehran, Iran.} \email{mtousi@ipm.ir}


\keywords {Koszul homology module; local cohomology module; local homology module; Serre class. \\
The research of the first author is supported by grants from IPM (No. 95130212).}

\begin{abstract}
Given a Serre class $\mathcal{S}$ of modules, we compare the containment of the Koszul homology, Ext modules, Tor modules, local homology, and local cohomology in $\mathcal{S}$ up to a given bound $s \geq 0$. As some applications, we give a full characterization of noetherian local homology modules. Further, we establish a comprehensive vanishing result which readily leads to the formerly known descriptions of the numerical invariants width and depth in terms of Koszul homology, local homology, and local cohomology. Also, we immediately recover a few renowned vanishing criteria scattered about the literature.
\end{abstract}

\maketitle

\tableofcontents

\section{Introduction}

\sloppy

Throughout this paper, $R$ denotes a commutative noetherian ring with identity, and $\mathcal{M}(R)$ designates the category of $R$-modules.

The theory of local cohomology was invented and introduced by A. Grothendieck about six decades ago and has since been tremendously delved and developed by various algebraic geometers and commutative algebraists.

There exists a dual theory, called the theory of local homology, which was initiated by Matlis \cite{Mat2} and \cite{Mat1} in 1974, and its investigation was continued by Simon in \cite{Si1} and \cite{Si2}. After the ground-breaking achievements of Greenlees and May \cite{GM}, and Alonso Tarr\'{i}o, Jerem\'{i}as L\'{o}pez and Lipman \cite {AJL}, a new era in the study of this theory has undoubtedly begun; see e.g. \cite{Sc}, \cite{CN1}, \cite{CN2},
\cite{Fr}, \cite{Ri}, \cite{MD}, \cite{HD} and \cite{DFT}. However, the theory of local homology has not been addressed as it merits and numerous counterparts of the results in local cohomology are yet to be explored and unearthed in local homology.

One interesting problem in local cohomology is the investigation of the artinian property. This is very well accomplished in \cite[Theorem 5.5]{Me2}. The dual problem, i.e. the noetherian property of local homology will be, among other things, probed carefully in this paper.

An excursion among the results \cite[Theorem 2.10]{AM1}, \cite[Theorem 2.9]{AM2}, \cite[Corollary 3.1]{BA}, \cite[Lemma 2.1]{BKN}, \cite[Proposition 1 and Corollary 1]{DM}, \cite[Lemma 4.2]{HK}, \cite[Theorem 2.1]{Me2}, \cite[Propositions 7.1, 7.2 and 7.4]{WW}, and \cite[Lemma 1.2]{Yo} reveals an in-depth connection between local cohomology, Ext modules, Tor modules, and Koszul homology in terms of their containment in a Serre class of modules. The purpose of this paper is to bring the local homology into play and uncover the true connection between all these homology and cohomology theories, and consequently, to illuminate and enhance the aforementioned results.

To shed some light on this revelation, we observe that given elements $\underline{a}= a_{1},\ldots,a_{n} \in R$ and an $R$-module $M$, the Koszul homology $H_{i}(\underline{a};M)=0$ for every $i < 0$ or $i > n$, and further
$$H_{n}(\underline{a};M)\cong \left(0:_{M}(\underline{a})\right) \cong \Ext_{R}^{0}\left(R/(\underline{a}),M\right)$$
and
$$H_{0}(\underline{a};M)\cong M/(\underline{a})M \cong \Tor^{R}_{0}\left(R/(\underline{a}),M\right).$$
These isomorphisms suggest the existence of an intimate connection between the Koszul homology on the one hand, and Ext and Tor modules on the other hand. The connection seems to manifest as in the following casual diagram for any given integer $s \geq 0$:\\

\begin{center}
\begin{tikzpicture}
  \node(1) at (-2.23,1.5) {$\underbrace{\Ext_{R}^{0}\left(R/(\underline{a}),M\right),\ldots,\Ext_{R}^{s}\left(R/(\underline{a}),M\right)}$};
  \node(2) at (-2.23,0.78) {$\updownarrow$};
  \node(3) at (0,0) {$\overbrace{H_{n}(\underline{a};M),\ldots,H_{n-s}(\underline{a};M)},\ldots,\underbrace{H_{s}(\underline{a};M),\ldots,H_{0}(\underline{a};M)}$};
  \node(4) at (2.43,-0.78) {$\updownarrow$};
  \node(5) at (2.43,-1.5) {$\overbrace{\Tor_{s}^{R}\left(R/(\underline{a}),M\right),\ldots,\Tor_{0}^{R}\left(R/(\underline{a}),M\right)}$};
\end{tikzpicture}
\end{center}

As the above diagram depicts, the Koszul homology acts as a connecting device in attaching the Ext modules to Tor modules in an elegant way. One such connection, one may surmise, could be the containment in some Serre class of modules that we are about to study in the following sections. More specifically, we prove the following as our main results; see Definitions \ref{2.1}, \ref{2.3}, \ref{2.5}, and \ref{2.7}.

\begin{theorem} \label{1.1}
Let $\mathfrak{a}=(a_{1},\ldots,a_{n})$ be an ideal of $R$, $\underline{a}=a_{1},\ldots,a_{n}$, and $M$ an $R$-module. Let $\mathcal{P}$ be a Serre property. Then the following three conditions are equivalent for any given $s \geq 0$:
\begin{enumerate}
\item[(i)] $H_{i}(\underline{a};M) \in \mathcal{S}_{\mathcal{P}}(R)$ for every $0 \leq i \leq s$.
\item[(ii)] $\Tor^{R}_{i}(N,M) \in \mathcal{S}_{\mathcal{P}}(R)$ for every finitely generated $R$-module $N$ with $\Supp_{R}(N)\subseteq V(\mathfrak{a})$, and for every $0 \leq i \leq s$.
\item[(iii)] $\Tor^{R}_{i}(N,M) \in \mathcal{S}_{\mathcal{P}}(R)$ for some finitely generated $R$-module $N$ with $\Supp_{R}(N)= V(\mathfrak{a})$, and for every $0 \leq i \leq s$.
\end{enumerate}
If in addition, $\mathcal{P}$ satisfies the condition $\mathfrak{D}_{\mathfrak{a}}$, then the above three conditions are equivalent to the following condition:
\begin{enumerate}
\item[(iv)] $H^{\mathfrak{a}}_{i}(M)\in \mathcal{S}_{\mathcal{P}}\left(\widehat{R}^{\mathfrak{a}}\right)$ for every $0 \leq i \leq s$.
\end{enumerate}
\end{theorem}

And dually:

\begin{theorem} \label{1.2}
Let $\mathfrak{a}=(a_{1},\ldots,a_{n})$ be an ideal of $R$, $\underline{a}=a_{1},\ldots,a_{n}$, and $M$ an $R$-module. Let $\mathcal{P}$ be a Serre property. Then the following three conditions are equivalent for any given $s \geq 0$:
\begin{enumerate}
\item[(i)] $H_{n-i}(\underline{a};M) \in \mathcal{S}_{\mathcal{P}}(R)$ for every $0 \leq i \leq s$.
\item[(ii)] $\Ext^{i}_{R}(N,M) \in \mathcal{S}_{\mathcal{P}}(R)$ for every finitely generated $R$-module $N$ with $\Supp_{R}(N)\subseteq V(\mathfrak{a})$, and for every $0 \leq i \leq s$.
\item[(iii)] $\Ext^{i}_{R}(N,M) \in \mathcal{S}_{\mathcal{P}}(R)$ for some finitely generated $R$-module $N$ with $\Supp_{R}(N)= V(\mathfrak{a})$, and for every $0 \leq i \leq s$.
\end{enumerate}
If in addition, $\mathcal{P}$ satisfies the condition $\mathfrak{C}_{\mathfrak{a}}$, then the above three conditions are equivalent to the following condition:
\begin{enumerate}
\item[(iv)] $H^{i}_{\mathfrak{a}}(M)\in \mathcal{S}_{\mathcal{P}}(R)$ for every $0 \leq i \leq s$.
\end{enumerate}
\end{theorem}

Having these results at our disposal, by specializing when $s = n$ and the Serre property is noetherianness, artinianness, or being zero, we draw a handful of fruitful conclusions on noetherianness, artinianness, and vanishing of the five foregoing homology and cohomology theories.

\section{Prerequisites}

First and foremost, we recall the notion of a Serre class of modules.

\begin{definition} \label{2.1}
A class $\mathcal{S}$ of $R$-modules is said to be a \textit{Serre class} if, given any short exact sequence
$$0 \rightarrow M^{\prime} \rightarrow M \rightarrow M^{\prime\prime} \rightarrow 0$$
of $R$-modules, we have $M \in \mathcal{S}$ if and only if $M^{\prime}, M^{\prime\prime} \in \mathcal{S}$.
\end{definition}

The following example showcases the stereotypical instances of Serre classes.

\begin{example} \label{2.2}
Let $\mathfrak{a}$ be an ideal of $R$. Then the following classes of modules are Serre classes:
\begin{enumerate}
\item[(i)] The zero class.
\item[(ii)] The class of all noetherian $R$-modules.
\item[(iii)] The class of all artinian $R$-modules.
\item[(iv)] The class of all minimax $R$-modules.
\item[(v)] The class of all minimax and $\mathfrak{a}$-cofinite $R$-modules; see \cite[Corollary 4.4]{Me1}.
\item[(vi)] The class of all weakly Laskerian $R$-modules.
\item[(v)] The class of all Matlis reflexive $R$-modules.
\end{enumerate}
\end{example}

We will be concerned with a change of rings when dealing with Serre classes. Therefore, we adopt the following notion of Serre property to exclude any possible ambiguity in the statement of the results.

\begin{definition} \label{2.3}
A property $\mathcal{P}$ concerning modules is said to be a \textit{Serre property} if
$$\mathcal{S}_{\mathcal{P}}(R):= \left\{M\in \mathcal{M}(R) \suchthat M \text{ satisfies the property } \mathcal{P} \right\}$$
is a Serre class for every ring $R$.
\end{definition}

Given a Serre property, there naturally arise the corresponding notions of depth and width.

\begin{definition} \label{2.4}
Let $\mathfrak{a}$ be an ideal of $R$, and $M$ an $R$-module. Let $\mathcal{P}$ be a Serre property. Then:
\begin{enumerate}
\item[(i)] Define the $\mathcal{P}$-depth of $M$ with respect to $\mathfrak{a}$ to be
\begin{center}
$\mathcal{P}$-$\depth_{R}(\mathfrak{a},M):= \inf \left\{i\geq 0 \suchthat \Ext^{i}_{R}(R/\mathfrak{a},M) \notin \mathcal{S}_{\mathcal{P}}(R) \right\}.$
\end{center}
\item[(ii)] Define the $\mathcal{P}$-width of $M$ with respect to $\mathfrak{a}$ to be
\begin{center}
$\mathcal{P}$-$\width_{R}(\mathfrak{a},M):= \inf \left\{i\geq 0 \suchthat \Tor^{R}_{i}(R/\mathfrak{a},M) \notin \mathcal{S}_{\mathcal{P}}(R) \right\}.$
\end{center}
\end{enumerate}
\end{definition}

It is clear that upon letting $\mathcal{P}$ be the Serre property of being zero, we recover the classical notions of depth and width.

We next remind the definition of local cohomology and local homology functors. Let $\mathfrak{a}$ be an ideal of $R$. We let
$$\Gamma_{\mathfrak{a}}(M):=\left\{x\in M \suchthat \mathfrak{a}^{t}x=0 \text{ for some } t\geq 1 \right\}$$
for an $R$-module $M$, and $\Gamma_{\mathfrak{a}}(f):=f|_{\Gamma_{\mathfrak{a}}(M)}$ for an $R$-homomorphism $f:M\rightarrow N$. This provides us with the so-called $\mathfrak{a}$-torsion functor $\Gamma_{\mathfrak{a}}(-)$ on the category of $R$-modules. The $i$th local cohomology functor with respect to $\mathfrak{a}$ is defined to be
$$H^{i}_{\mathfrak{a}}(-):= R^{i}\Gamma_{\mathfrak{a}}(-)$$
for every $i \geq 0$. In addition, the cohomological dimension of $\mathfrak{a}$ with respect to $M$ is
$$\cd(\mathfrak{a},M) := \sup \left\{i \geq 0 \suchthat H^{i}_{\mathfrak{a}}(M)\neq 0 \right\}.$$

Recall that an $R$-module $M$ is said to be $\mathfrak{a}$-torsion if $M=\Gamma_{\mathfrak{a}}(M)$. It is well-known that any $\mathfrak{a}$-torsion $R$-module $M$ enjoys a natural $\widehat{R}^{\mathfrak{a}}$-module structure in such a way that the lattices of its $R$-submodules and $\widehat{R}^{\mathfrak{a}}$-submodules coincide. Besides, we have $\widehat{R}^{\mathfrak{a}}\otimes_{R}M \cong M$ both as $R$-modules and $\widehat{R}^{\mathfrak{a}}$-modules.

Likewise, we let
$$\Lambda^{\mathfrak{a}}(M):=\widehat{M}^{\mathfrak{a}}= \underset{n}\varprojlim (M/\mathfrak{a}^{n}M)$$
for an $R$-module $M$, and $\Lambda^{\mathfrak{a}}(f):=\widehat{f}$ for an $R$-homomorphism $f:M\rightarrow N$. This provides us with the so-called $\mathfrak{a}$-adic completion functor $\Lambda^{\mathfrak{a}}(-)$ on the category of $R$-modules. The $i$th local homology functor with respect to $\mathfrak{a}$ is defined to be
$$H^{\mathfrak{a}}_{i}(-):= L_{i}\Lambda^{\mathfrak{a}}(-)$$
for every $i \geq 0$. Besides, the homological dimension of $\mathfrak{a}$ with respect to $M$ is
$$\hd(\mathfrak{a},M) := \sup \left\{i \geq 0 \suchthat H^{\mathfrak{a}}_{i}(M)\neq 0 \right\}.$$

The next conditions are also required to be imposed on the Serre classes when we intend to bring the local homology and local cohomology into the picture.

\begin{definition} \label{2.5}
Let $\mathcal{P}$ be a Serre property, and $\mathfrak{a}$ an ideal of $R$. We say that $\mathcal{P}$ satisfies the condition $\mathfrak{D}_{\mathfrak{a}}$ if the following statements hold:
\begin{enumerate}
\item[(i)] If $R$ is $\mathfrak{a}$-adically complete and $M/\mathfrak{a}M \in \mathcal{S}_{\mathcal{P}}(R)$ for some $R$-module $M$, then $H^{\mathfrak{a}}_{0}(M) \in \mathcal{S}_{\mathcal{P}}(R)$.
\item[(ii)] For any $\mathfrak{a}$-torsion $R$-module $M$, we have $M\in \mathcal{S}_{\mathcal{P}}(R)$ if and only if $M\in \mathcal{S}_{\mathcal{P}}\left(\widehat{R}^{\mathfrak{a}}\right)$.
\end{enumerate}
\end{definition}

\begin{example} \label{2.6}
Let $\mathfrak{a}$ be an ideal of $R$. Then we have:
\begin{enumerate}
\item[(i)] The Serre property of being zero satisfies the condition $\mathfrak{D}_{\mathfrak{a}}$. Use \cite[Lemma 5.1 (iii)]{Si1}.
\item[(ii)] The Serre property of being noetherian satisfies the condition $\mathfrak{D}_{\mathfrak{a}}$. See \cite[Lemma 2.5]{DFT}.
\end{enumerate}
\end{example}

Melkersson \cite[Definition 2.1]{AM3} defines the condition  $\mathfrak{C}_{\mathfrak{a}}$  for Serre classes. Adopting his definition, we have:

\begin{definition} \label{2.7}
Let $\mathcal{P}$ be a Serre property, and $\mathfrak{a}$ an ideal of $R$. We say that $\mathcal{P}$ satisfies the condition $\mathfrak{C}_{\mathfrak{a}}$ if the containment $(0:_{M}\mathfrak{a}) \in \mathcal{S}_{\mathcal{P}}(R)$ for some $R$-module $M$, implies that $\Gamma_{\mathfrak{a}}(M) \in \mathcal{S}_{\mathcal{P}}(R)$.
\end{definition}

\begin{example} \label{2.8}
Let $\mathfrak{a}$ be an ideal of $R$. Then we have:
\begin{enumerate}
\item[(i)] The Serre property of being zero satisfies the condition $\mathfrak{C}_{\mathfrak{a}}$. Indeed, if $M$ is an $R$-module such that $(0:_{M}\mathfrak{a})=0$, then it can be easily seen that $\Gamma_{\mathfrak{a}}(M)=0$.
\item[(ii)] The Serre property of being artinian satisfies the condition $\mathfrak{C}_{\mathfrak{a}}$. This follows from the Melkersson's Criterion \cite[Theorem 1.3]{Me2}.
\item[(iii)] A Serre class which is closed under taking direct limits is called a torsion theory. There exists a whole bunch of torsion theories, e.g. given any $R$-module $L$, the class
    $$\mathfrak{T}_{L}:= \left\{M\in \mathcal{M}(R) \suchthat \Supp_{R}(M)\subseteq \Supp_{R}(L) \right\}$$
    is a torsion theory. One can easily see that any torsion theory satisfies the condition $\mathfrak{C}_{\mathfrak{a}}$.
\end{enumerate}
\end{example}

We finally describe the Koszul homology briefly. The Koszul complex $K^{R}(a)$ on an element $a \in R$ is the $R$-complex
$$K^{R}(a):=\Cone(R \xrightarrow{a} R),$$
and the Koszul complex $K^{R}(\underline{a})$ on a sequence of elements $\underline{a} = a_{1},\ldots,a_{n} \in R$ is the $R$-complex
$$K^{R}(\underline{a}):= K^{R}(a_{1}) \otimes_{R} \cdots \otimes_{R} K^{R}(a_{n}).$$
It is easy to see that $K^{R}(\underline{a})$ is a complex of finitely generated free $R$-modules concentrated in degrees $n,\ldots,0$.
Given any $R$-module $M$, there is an isomorphism of $R$-complexes
$$K^{R}(\underline{a})\otimes_{R}M \cong \Sigma^{n} \Hom_{R}\left(K^{R}(\underline{a}),M\right),$$
which is sometimes referred to as the self-duality property of the Koszul complex. Accordingly, we feel free to define the Koszul homology of the sequence $\underline{a}$ with coefficients in $M$, by setting
$$H_{i}(\underline{a};M):= H_{i}\left(K^{R}(\underline{a})\otimes_{R}M\right) \cong H_{i-n}\left(\Hom_{R}\left(K^{R}(\underline{a}),M\right)\right)$$
for every $i\geq 0$.

\section{Proofs of the main results}

In this section, we study the containment of Koszul homology, Ext modules, Tor modules, local homology, and local cohomology in Serre classes.

In the proof of Theorem \ref{1.1} below, we use the straightforward observation that given elements $\underline{a}= a_{1},\dots,a_{n} \in R$, a finitely generated $R$-module $N$, and an $R$-module $M$, if $M$ belongs to a Serre class $\mathcal{S}$, then $H_{i}(\underline{a};M)\in \mathcal{S}$, $\Ext^{i}_{R}(N,M) \in \mathcal{S}$, and $\Tor^{R}_{i}(N,M) \in \mathcal{S}$ for every $i \geq 0$. For a proof refer to \cite[Lemma 2.1]{AT}.

{\bf Proof of Theorem 1.1.}
(i) $\Rightarrow$ (iii): Set $N=R/\mathfrak{a}$. Let $F$ be a free resolution of $R/\mathfrak{a}$ consisting of finitely generated $R$-modules. Then the $R$-complex $F\otimes_{R}M$ is isomorphic to an $R$-complex of the form
$$\cdots \rightarrow M^{s_{2}} \xrightarrow{\partial_{2}} M^{s_{1}} \xrightarrow{\partial_{1}} M^{s_{0}} \rightarrow 0.$$
We note that
\begin{equation} \label{eq:3.1.1}
\Tor_{0}^{R}(R/\mathfrak{a},M) \cong \coker \partial_{1} \cong H_{0}(\underline{a};M) \in \mathcal{S}_{\mathcal{P}}(R),
\end{equation}
by the assumption. If $s=0$, then we are done. Suppose that $s \geq 1$. The short exact sequence
$$0 \rightarrow \im \partial_{1} \rightarrow M^{s_{0}} \rightarrow \coker \partial_{1} \rightarrow 0,$$
induces the exact sequence
$$H_{i+1}\left(\underline{a};\coker \partial_{1}\right) \rightarrow H_{i}\left(\underline{a};\im \partial_{1}\right) \rightarrow H_{i}(\underline{a};M^{s_{0}}).$$
The assumption together with the display \eqref{eq:3.1.1} imply that the two lateral terms of the above exact sequence belong to $\mathcal{S}_{\mathcal{P}}(R)$, so $H_{i}\left(\underline{a};\im \partial_{1}\right) \in \mathcal{S}_{\mathcal{P}}(R)$ for every $0 \leq i \leq s$. The short exact sequence
$$0 \rightarrow \ker \partial_{1} \rightarrow M^{s_{1}} \rightarrow \im \partial_{1} \rightarrow 0,$$
yields the exact sequence
$$H_{i+1}\left(\underline{a};\im \partial_{1}\right) \rightarrow H_{i}\left(\underline{a};\ker \partial_{1}\right) \rightarrow H_{i}(\underline{a};M^{s_{1}}).$$
Therefore, $H_{i}\left(\underline{a};\ker \partial_{1}\right) \in \mathcal{S}_{\mathcal{P}}(R)$ for every $0 \leq i \leq s-1$. As $s \geq 1$, we see that $H_{0}\left(\underline{a};\ker \partial_{1}\right) \in \mathcal{S}_{\mathcal{P}}(R)$. The short exact sequence
\begin{equation} \label{eq:3.1.2}
0 \rightarrow \im \partial_{2} \rightarrow \ker \partial_{1} \rightarrow \Tor_{1}^{R}(R/\mathfrak{a},M) \rightarrow 0,
\end{equation}
implies the exact sequence
$$H_{0}\left(\underline{a};\ker \partial_{1}\right) \rightarrow H_{0}\left(\underline{a};\Tor_{1}^{R}(R/\mathfrak{a},M)\right) \rightarrow 0.$$
Therefore, $H_{0}\left(\underline{a};\Tor_{1}^{R}(R/\mathfrak{a},M)\right) \in \mathcal{S}_{\mathcal{P}}(R)$.
But $\mathfrak{a} \Tor^{R}_{1}(R/\mathfrak{a},M)=0$, so
$$\Tor^{R}_{1}(R/\mathfrak{a},M) \cong H_{0}\left(\underline{a};\Tor_{1}^{R}(R/\mathfrak{a},M)\right) \in \mathcal{S}_{\mathcal{P}}(R).$$
If $s=1$, then we are done. Suppose that $s \geq 2$. The short exact sequence \eqref{eq:3.1.2} induces the exact sequence
$$H_{i+1}\left(\underline{a};\Tor_{1}^{R}(R/\mathfrak{a},M)\right) \rightarrow H_{i}\left(\underline{a};\im \partial_{2}\right) \rightarrow H_{i}\left(\underline{a};\ker \partial_{1}\right).$$
It follows that $H_{i}\left(\underline{a};\im \partial_{2}\right) \in \mathcal{S}_{\mathcal{P}}(R)$ for every $0 \leq i \leq s-1$. The short exact sequence
$$0 \rightarrow \ker \partial_{2} \rightarrow M^{s_{2}} \rightarrow \im \partial_{2} \rightarrow 0,$$
yields the exact sequence
$$H_{i+1}\left(\underline{a};\im \partial_{2}\right) \rightarrow H_{i}\left(\underline{a};\ker \partial_{2}\right) \rightarrow H_{i}(\underline{a};M^{s_{2}}).$$
Thus $H_{i}\left(\underline{a};\ker \partial_{2}\right) \in \mathcal{S}_{\mathcal{P}}(R)$ for every $0 \leq i \leq s-2$.
As $s \geq 2$, we see that $H_{0}\left(\underline{a};\ker \partial_{2}\right) \in \mathcal{S}_{\mathcal{P}}(R)$. The short exact sequence
$$0 \rightarrow \im \partial_{3} \rightarrow \ker \partial_{2} \rightarrow \Tor_{2}^{R}\left(R/\mathfrak{a},M\right) \rightarrow 0,$$
yields the exact sequence
$$H_{0}\left(\underline{a};\ker \partial_{2}\right) \rightarrow H_{0}\left(\underline{a};\Tor_{2}^{R}(R/\mathfrak{a},M)\right) \rightarrow 0.$$
As before, we conclude that
$$\Tor^{R}_{2}(R/\mathfrak{a},M) \cong H_{0}\left(\underline{a};\Tor_{2}^{R}(R/\mathfrak{a},M)\right) \in \mathcal{S}_{\mathcal{P}}(R).$$
If $s=2$, then we are done. Proceeding in this manner, we see that $\Tor^{R}_{i}(R/\mathfrak{a},M) \in \mathcal{S}_{\mathcal{P}}(R)$ for every $0 \leq i \leq s$.

(iii) $\Rightarrow$ (ii): Let $L$ be a finitely generated $R$-module with $\Supp_{R}(L) \subseteq V(\mathfrak{a})$. By induction on $s$, we show that $\Tor_{i}^{R}(L,M)\in \mathcal{S}_{\mathcal{P}}(R)$ for every $0 \leq i \leq s$. Let $s=0$. Then $N \otimes_{R} M \in \mathcal{S}(R)$. Using Gruson's Filtration Lemma \cite[Theorem 4.1]{Va}, there is a finite filtration
$$0=L_{0} \subseteq L_{1} \subseteq \cdots \subseteq L_{t-1} \subseteq L_{t} = L,$$
such that $L_{j}/L_{j-1}$ is isomorphic to a quotient of a finite direct sum of copies of $N$ for every $1 \leq j \leq t$. Given any $1 \leq j \leq t$, the exact sequence
\begin{equation} \label{eq:3.1.3}
0 \rightarrow K_{j} \rightarrow N^{r_{j}} \rightarrow L_{j}/L_{j-1} \rightarrow 0,
\end{equation}
induces the exact sequence
$$N^{r_{j}} \otimes_{R} M \rightarrow (L_{j}/L_{j-1}) \otimes_{R} M \rightarrow 0.$$
It follows that $(L_{j}/L_{j-1}) \otimes_{R} M \in \mathcal{S}_{\mathcal{P}}(R)$. Now the short exact sequence
\begin{equation} \label{eq:3.1.4}
0 \rightarrow L_{j-1} \rightarrow L_{j} \rightarrow L_{j}/L_{j-1} \rightarrow 0,
\end{equation}
yields the exact sequence
$$L_{j-1} \otimes_{R} M \rightarrow L_{j} \otimes_{R} M \rightarrow (L_{j}/L_{j-1}) \otimes_{R} M.$$
A successive use of the above exact sequence, letting $j=1,\ldots,t$, implies that $L\otimes_{R}M \in \mathcal{S}_{\mathcal{P}}(R)$.

Now let $s \geq 1$, and suppose that the results holds true for $s-1$. The induction hypothesis implies that $\Tor^{R}_{i}(L,M) \in \mathcal{S}_{\mathcal{P}}(R)$ for every $0 \leq i \leq s-1$. Hence, it suffices to show that $\Tor^{R}_{s}(L,M) \in \mathcal{S}_{\mathcal{P}}(R)$. Given any $1 \leq j \leq t$, the short exact sequence \eqref{eq:3.1.3} induces the exact sequence $$\Tor^{R}_{s}(N^{r_{j}},M) \rightarrow \Tor^{R}_{s}\left(L_{j}/L_{j-1},M\right) \rightarrow \Tor^{R}_{s-1}(K_{j},M).$$
The induction hypothesis shows that $\Tor^{R}_{s-1}(K_{j},M) \in \mathcal{S}_{\mathcal{P}}(R)$, so from the above exact sequence we get that $\Tor^{R}_{s}\left(L_{j}/L_{j-1},M\right) \in \mathcal{S}_{\mathcal{P}}(R)$. Now, the short exact sequence \eqref{eq:3.1.4} yields the exact sequence
$$\Tor^{R}_{s}(L_{j-1},M) \rightarrow \Tor^{R}_{s}(L_{j},M) \rightarrow \Tor^{R}_{s}\left(L_{j}/L_{j-1},M\right).$$
A successive use of the above exact sequence, letting $j=1,\ldots,t$, implies that $\Tor^{R}_{s}(L,M) \in \mathcal{S}_{\mathcal{P}}(R)$.

(ii) $\Rightarrow$ (i): The hypothesis implies that $\Tor^{R}_{i}(R/\mathfrak{a},M) \in \mathcal{S}_{\mathcal{P}}(R)$ for every $0 \leq i \leq s$. The $R$-complex $K^{R}(\underline{a})\otimes_{R}M$ is isomorphic to an $R$-complex of the form
$$0 \rightarrow M \xrightarrow{\partial_{n}} M^{t_{n-1}} \rightarrow \cdots \rightarrow M^{t_{2}} \xrightarrow{\partial_{2}} M^{t_{1}} \xrightarrow{\partial_{1}} M \rightarrow 0.$$
We have
\begin{equation} \label{eq:3.1.5}
H_{0}(\underline{a};M)\cong \coker \partial_{1} \cong \Tor^{R}_{0}(R/\mathfrak{a},M) \in \mathcal{S}_{\mathcal{P}}(R),
\end{equation}
by the assumption. If $s=0$, then we are done. Suppose that $s \geq 1$. The short exact sequence
$$0 \rightarrow \im \partial_{1} \rightarrow M \rightarrow \coker \partial_{1} \rightarrow 0,$$
induces the exact sequence
$$\Tor^{R}_{i+1}(R/\mathfrak{a},\coker \partial_{1}) \rightarrow \Tor^{R}_{i}(R/\mathfrak{a},\im \partial_{1}) \rightarrow \Tor^{R}_{i}(R/\mathfrak{a},M).$$
The assumption along with the display \eqref{eq:3.1.5} imply that the two lateral terms of the above exact sequence belong to $\mathcal{S}_{\mathcal{P}}(R)$, so $\Tor^{R}_{i}(R/\mathfrak{a},\im \partial_{1}) \in \mathcal{S}_{\mathcal{P}}(R)$ for every $0 \leq i \leq s$. The short exact sequence
$$0 \rightarrow \ker \partial_{1} \rightarrow M^{t_{1}} \rightarrow \im \partial_{1} \rightarrow 0,$$
yields the exact sequence
$$\Tor^{R}_{i+1}(R/\mathfrak{a},\im \partial_{1}) \rightarrow \Tor^{R}_{i}(R/\mathfrak{a},\ker \partial_{1}) \rightarrow \Tor^{R}_{i}(R/\mathfrak{a},M^{t_{1}}).$$
Therefore, $\Tor^{R}_{i}(R/\mathfrak{a},\ker \partial_{1}) \in \mathcal{S}_{\mathcal{P}}(R)$ for every $0 \leq i \leq s-1$.
As $s \geq 1$, we see that $(R/\mathfrak{a})\otimes_{R} \ker \partial_{1} \in \mathcal{S}_{\mathcal{P}}(R)$. The short exact sequence
\begin{equation} \label{eq:3.1.6}
0 \rightarrow \im \partial_{2} \rightarrow \ker \partial_{1} \rightarrow H_{1}(\underline{a};M) \rightarrow 0,
\end{equation}
implies the exact sequence
$$(R/\mathfrak{a})\otimes_{R} \ker \partial_{1} \rightarrow (R/\mathfrak{a})\otimes_{R} H_{1}(\underline{a};M) \rightarrow 0.$$
Therefore, $(R/\mathfrak{a})\otimes_{R} H_{1}(\underline{a};M) \in \mathcal{S}_{\mathcal{P}}(R)$.
But $\mathfrak{a} H_{1}(\underline{a};M)=0$, so
$$H_{1}(\underline{a};M) \cong (R/\mathfrak{a})\otimes_{R} H_{1}(\underline{a};M) \in \mathcal{S}_{\mathcal{P}}(R).$$
If $s=1$, then we are done. Suppose that $s \geq 2$. The short exact sequence \eqref{eq:3.1.6} induces the exact sequence
$$\Tor^{R}_{i+1}\left(R/\mathfrak{a},H_{1}(\underline{a};M)\right) \rightarrow \Tor^{R}_{i}(R/\mathfrak{a},\im \partial_{2}) \rightarrow \Tor^{R}_{i}(R/\mathfrak{a},\ker \partial_{1}).$$
It follows that $\Tor^{R}_{i}(R/\mathfrak{a},\im \partial_{2}) \in \mathcal{S}_{\mathcal{P}}(R)$ for every $0 \leq i \leq s-1$. The short exact sequence
$$0 \rightarrow \ker \partial_{2} \rightarrow M^{t_{2}} \rightarrow \im \partial_{2} \rightarrow 0,$$
yields the exact sequence
$$\Tor^{R}_{i+1}\left(R/\mathfrak{a},\im \partial_{2}\right) \rightarrow \Tor^{R}_{i}(R/\mathfrak{a},\ker \partial_{2}) \rightarrow \Tor^{R}_{i}(R/\mathfrak{a},M^{t_{2}}).$$
Thus $\Tor^{R}_{i}(R/\mathfrak{a},\ker \partial_{2}) \in \mathcal{S}_{\mathcal{P}}(R)$ for every $0 \leq i \leq s-2$. As $s \geq 2$, we see that $(R/\mathfrak{a}) \otimes_{R} \ker \partial_{2} \in \mathcal{S}_{\mathcal{P}}(R)$. The short exact sequence
$$0 \rightarrow \im \partial_{3} \rightarrow \ker \partial_{2} \rightarrow H_{2}(\underline{a};M) \rightarrow 0,$$
implies the exact sequence
$$(R/\mathfrak{a}) \otimes_{R} \ker \partial_{2} \rightarrow (R/\mathfrak{a}) \otimes_{R} H_{2}(\underline{a};M) \rightarrow 0.$$
As before, we conclude that
$$H_{2}(\underline{a};M) \cong (R/\mathfrak{a}) \otimes_{R} H_{2}(\underline{a};M) \in \mathcal{S}_{\mathcal{P}}(R).$$
If $s=2$, then we are done. Proceeding in this manner, we infer that $H_{i}(\underline{a};M) \in \mathcal{S}_{\mathcal{P}}(R)$ for every $0 \leq i \leq s$.

Now, assume that the Serre property $\mathcal{P}$ satisfies the condition $\mathfrak{D}_{\mathfrak{a}}$. We first note that since $\Tor^{R}_{i}(R/\mathfrak{a},M)$ is an $\mathfrak{a}$-torsion $R$-module, it has an $\widehat{R}^{\mathfrak{a}}$-module structure such that
$$\Tor^{R}_{i}(R/\mathfrak{a},M) \cong \widehat{R}^{\mathfrak{a}} \otimes_{R} \Tor^{R}_{i}(R/\mathfrak{a},M) \cong \Tor^{\widehat{R}^{\mathfrak{a}}}_{i}\left(\widehat{R}^{\mathfrak{a}}/\mathfrak{a}\widehat{R}^{\mathfrak{a}},\widehat{R}^{\mathfrak{a}} \otimes_{R} M\right),$$
for every $i \geq 0$ both as $R$-modules and $\widehat{R}^{\mathfrak{a}}$-modules. Moreover, by \cite[Lemma 2.3]{Si2}, we have
$$H^{\mathfrak{a}}_{i}(M)\cong H^{\mathfrak{a}\widehat{R}^{\mathfrak{a}}}_{i}\left(\widehat{R}^{\mathfrak{a}}\otimes_{R}M\right),$$
for every $i \geq 0$ both as $R$-modules and $\widehat{R}^{\mathfrak{a}}$-modules. With this preparation we prove:

(iv) $\Rightarrow$ (iii): We let $N=R/\mathfrak{a}$, and show that $\Tor_{i}^{R}(R/\mathfrak{a},M)\in \mathcal{S}_{\mathcal{P}}(R)$ for every $0 \leq i \leq s$. First suppose that $R$ is $\mathfrak{a}$-adically complete. By \cite[Lemma 2.6]{DFT}, there is a first quadrant spectral sequence
\begin{equation} \label{eq:3.1.7}
E^{2}_{p,q}= \Tor^{R}_{p}\left(R/\mathfrak{a},H^{\mathfrak{a}}_{q}(M)\right) \underset {p} \Rightarrow \Tor^{R}_{p+q}(R/\mathfrak{a},M).
\end{equation}
The hypothesis implies that $E^{2}_{p,q}\in \mathcal{S}_{\mathcal{P}}(R)$ for every $p \geq 0$ and $0 \leq q \leq s$. Let $0 \leq i \leq s$. There is a finite filtration
$$0=U^{-1} \subseteq U^{0} \subseteq \cdots \subseteq U^{i} = \Tor^{R}_{i}(R/\mathfrak{a},M),$$
such that $U^{p}/U^{p-1}\cong E^{\infty}_{p,i-p}$ for every $0 \leq p \leq i$. Since $E^{\infty}_{p,i-p}$ is a subquotient of $E^{2}_{p,i-p}$ and $0 \leq i-p \leq s$, we infer that
$$U^{p}/U^{p-1}\cong E^{\infty}_{p,i-p}\in \mathcal{S}_{\mathcal{P}}(R)$$
for every $0 \leq p \leq i$.
A successive use of the short exact sequence
$$0 \rightarrow U^{p-1} \rightarrow U^{p} \rightarrow U^{p}/U^{p-1} \rightarrow 0,$$
by letting $p=0,\ldots,i$, implies that $\Tor^{R}_{i}(R/\mathfrak{a},M)\in \mathcal{S}_{\mathcal{P}}(R)$.

Now, consider the general case. Since
$$H^{\mathfrak{a}\widehat{R}^{\mathfrak{a}}}_{i}\left(\widehat{R}^{\mathfrak{a}}\otimes_{R}M\right) \cong H^{\mathfrak{a}}_{i}(M) \in \mathcal{S}_{\mathcal{P}}\left(\widehat{R}^{\mathfrak{a}}\right),$$
for every $0 \leq i \leq s$, the special case implies that
$$\Tor^{R}_{i}(R/\mathfrak{a},M) \cong \Tor^{\widehat{R}^{\mathfrak{a}}}_{i}\left(\widehat{R}^{\mathfrak{a}}/\mathfrak{a}\widehat{R}^{\mathfrak{a}},\widehat{R}^{\mathfrak{a}} \otimes_{R} M\right)\in \mathcal{S}_{\mathcal{P}}\left(\widehat{R}^{\mathfrak{a}}\right),$$
for every $0 \leq i \leq s$. However, $\Tor^{R}_{i}(R/\mathfrak{a},M)$ is $\mathfrak{a}$-torsion, so the condition $\mathfrak{D}_{\mathfrak{a}}$ implies that $\Tor^{R}_{i}(R/\mathfrak{a},M)\in \mathcal{S}_{\mathcal{P}}(R)$ for every $0 \leq i \leq s$.

(ii) $\Rightarrow$ (iv): It follows from the hypothesis that $\Tor^{R}_{i}(R/\mathfrak{a},M)\in \mathcal{S}_{\mathcal{P}}(R)$ for every $0 \leq i \leq s$. First suppose that $R$ is $\mathfrak{a}$-adically complete. We argue by induction on $s$. Let $s=0$. Then $M/\mathfrak{a}M \in \mathcal{S}_{\mathcal{P}}(R)$, whence $H^{\mathfrak{a}}_{0}(M)\in \mathcal{S}_{\mathcal{P}}(R)$ by the condition $\mathfrak{D}_{\mathfrak{a}}$.

Now let $s \geq 1$, and suppose that the result holds true for $s-1$. The induction hypothesis implies that $H^{\mathfrak{a}}_{i}(M)\in \mathcal{S}_{\mathcal{P}}(R)$ for every $0 \leq i \leq s-1$. Hence, it suffices to show that $H^{\mathfrak{a}}_{s}(M)\in \mathcal{S}_{\mathcal{P}}(R)$.
Consider the spectral sequence \eqref{eq:3.1.7}. By the hypothesis, $E^{2}_{p,q}\in \mathcal{S}_{\mathcal{P}}(R)$ for every $p \in \mathbb{Z}$ and $0 \leq q \leq s-1$.
There is a finite filtration
$$0=U^{-1} \subseteq U^{0} \subseteq \cdots \subseteq U^{s} = \Tor^{R}_{s}(R/\mathfrak{a},M),$$
such that $U^{p}/U^{p-1}\cong E^{\infty}_{p,s-p}$ for every $0 \leq p \leq s$. As $\Tor^{R}_{s}(R/\mathfrak{a},M)\in \mathcal{S}_{\mathcal{P}}(R)$,
we conclude that
$$E^{\infty}_{0,s}\cong U^{0}/U^{-1} \cong U^{0}\in \mathcal{S}_{\mathcal{P}}(R).$$
Let $r \geq 2$, and consider the differentials
$$E^{r}_{r,s-r+1}\xrightarrow {d^{r}_{r,s-r+1}} E^{r}_{0,s}\xrightarrow {d^{r}_{0,s}} E^{r}_{-r,s+r-1}=0.$$
Since $s-r+1 \leq s-1$ and $E^{r}_{r,s-r+1}$ is a subquotient of $E^{2}_{r,s-r+1}$, the hypothesis implies that $E^{r}_{r,s-r+1}\in \mathcal{S}_{\mathcal{P}}(R)$, and consequently $\im d^{r}_{r,s-r+1}\in \mathcal{S}_{\mathcal{P}}(R)$ for every $r\geq 2$.
We thus obtain
$$E^{r+1}_{0,s} \cong \ker d^{r}_{0,s}/ \im d^{r}_{r,s-r+1} = E^{r}_{0,s}/ \im d^{r}_{r,s-r+1},$$
and consequently a short exact sequence
\begin{equation} \label{eq:3.1.8}
0 \rightarrow \im d^{r}_{r,s-r+1} \rightarrow E^{r}_{0,s} \rightarrow E^{r+1}_{0,s} \rightarrow 0.
\end{equation}
There is an integer $r_{0}\geq 2$, such that $E^{\infty}_{0,s} \cong E^{r+1}_{0,s}$ for every $r\geq r_{0}$. It follows that $E^{r_{0}+1}_{0,s}\in \mathcal{S}_{\mathcal{P}}(R)$.
Now, the short exact sequence \eqref{eq:3.1.8} implies that $E^{r_{0}}_{0,s}\in \mathcal{S}_{\mathcal{P}}(R)$. Using the short exact sequence
\eqref{eq:3.1.8} inductively, we conclude that
$$H^{\mathfrak{a}}_{s}(M)/\mathfrak{a}H^{\mathfrak{a}}_{s}(M)\cong E^{2}_{0,s}\in \mathcal{S}_{\mathcal{P}}(R).$$
Therefore, by \cite[Lemma 2.4]{DFT} and the condition $\mathfrak{D}_{\mathfrak{a}}$, we get
$$H^{\mathfrak{a}}_{s}(M) \cong H^{\mathfrak{a}}_{0}\left(H^{\mathfrak{a}}_{s}(M)\right)\in \mathcal{S}_{\mathcal{P}}(R).$$

Now, consider the general case. Since $\Tor^{R}_{i}(R/\mathfrak{a},M)$ is an $\mathfrak{a}$-torsion $R$-module such that $\Tor^{R}_{i}(R/\mathfrak{a},M)\in \mathcal{S}_{\mathcal{P}}(R)$ for every $0 \leq i \leq s$, we deduce that
$$\Tor^{\widehat{R}^{\mathfrak{a}}}_{i}\left(\widehat{R}^{\mathfrak{a}}/\mathfrak{a}\widehat{R}^{\mathfrak{a}},\widehat{R}^{\mathfrak{a}} \otimes_{R} M\right) \cong \Tor^{R}_{i}(R/\mathfrak{a},M) \in \mathcal{S}_{\mathcal{P}}\left(\widehat{R}^{\mathfrak{a}}\right),$$
for every $0 \leq i \leq s$. But the special case yields that
$$H^{\mathfrak{a}}_{i}(M) \cong H^{\mathfrak{a}\widehat{R}^{\mathfrak{a}}}_{i}\left(\widehat{R}^{\mathfrak{a}}\otimes_{R}M\right) \in \mathcal{S}_{\mathcal{P}}\left(\widehat{R}^{\mathfrak{a}}\right),$$
for every $0 \leq i \leq s$.
$\Box$

The following special case may be of independent interest.

\begin{corollary} \label{3.2}
Let $\mathfrak{a}=(a_{1},\ldots,a_{n})$ be an ideal of $R$, $\underline{a}=a_{1},\ldots,a_{n}$, and $M$ an $R$-module. Let $\mathcal{P}$ be a Serre property satisfying the condition $\mathfrak{D}_{\mathfrak{a}}$. Then the following conditions are equivalent:
\begin{enumerate}
\item[(i)] $M/\mathfrak{a}M \in \mathcal{S}_{\mathcal{P}}(R)$.
\item[(ii)] $H^{\mathfrak{a}}_{0}(M) \in \mathcal{S}_{\mathcal{P}}\left(\widehat{R}^{\mathfrak{a}}\right)$.
\item[(iii)] $\widehat{M}^{\mathfrak{a}}\in \mathcal{S}_{\mathcal{P}}\left(\widehat{R}^{\mathfrak{a}}\right)$.
\end{enumerate}
\end{corollary}

\begin{prf}
(i) $\Leftrightarrow$ (ii): Follows from Theorem \ref{1.1} upon letting $s=0$.

(ii) $\Rightarrow$ (iii): By \cite[Lemma 5.1 (i)]{Si1}, the natural homomorphism $H^{\mathfrak{a}}_{0}(M) \rightarrow \widehat{M}^{\mathfrak{a}}$ is surjective. Thus the result follows.

(iii) $\Rightarrow$ (i): Since $\widehat{M}^{\mathfrak{a}}\in \mathcal{S}_{\mathcal{P}}\left(\widehat{R}^{\mathfrak{a}}\right)$, we see that $\widehat{M}^{\mathfrak{a}}/\mathfrak{a}\widehat{M}^{\mathfrak{a}} \in \mathcal{S}_{\mathcal{P}}\left(\widehat{R}^{\mathfrak{a}}\right)$. However, by \cite[Theorem 1.1]{Si1}, we have $\widehat{M}^{\mathfrak{a}}/\mathfrak{a}\widehat{M}^{\mathfrak{a}} \cong M/\mathfrak{a}M$. It follows that $M/\mathfrak{a}M \in \mathcal{S}_{\mathcal{P}}\left(\widehat{R}^{\mathfrak{a}}\right)$. But $M/\mathfrak{a}M$ is $\mathfrak{a}$-torsion, so by the condition $\mathfrak{D}_{\mathfrak{a}}$, we have $M/\mathfrak{a}M \in \mathcal{S}_{\mathcal{P}}(R)$.
\end{prf}

Using similar arguments, we can prove the dual result to Theorem \ref{1.1}. It is worth noting that Theorem
\ref{1.2} is proved in \cite[Theorem 2.9]{AM3}. However, the condition $\mathfrak{C}_{\mathfrak{a}}$ is assumed to be satisfied in all four statements, but here we only require that the condition $\mathfrak{C}_{\mathfrak{a}}$ is satisfied for the last statement. Moreover, the techniques used there are quite different than those used here.

{\bf Proof of Theorem 1.2.}
The proof is similar to that of Theorem \ref{1.1}. However, instead of the spectral sequence \eqref{eq:3.1.7}, we use the third quadrant spectral sequence
$$E^{2}_{p,q}= \Ext^{-p}_{R}\left(R/\mathfrak{a},H^{-q}_{\mathfrak{a}}(M)\right) \underset {p} \Rightarrow \Ext^{-p-q}_{R}(R/\mathfrak{a},M).$$
$\Box$

If we let the integer $s$ exhaust the whole nonzero range of Koszul homology, i.e. $s=n$, then we can effectively combine Theorems \ref{1.1} and \ref{1.2} to obtain the following result which in turn generalizes \cite[Corollary 3.1]{BA}.

\begin{corollary} \label{3.4}
Let $\mathfrak{a}=(a_{1},\ldots,a_{n})$ be an ideal of $R$, $\underline{a}=a_{1},\ldots,a_{n}$, and $M$ an $R$-module. Let $\mathcal{P}$ be a Serre property. Then the following conditions are equivalent:
\begin{enumerate}
\item[(i)] $H_{i}(\underline{a};M)\in \mathcal{S}_{\mathcal{P}}(R)$ for every $0 \leq i \leq n$.
\item[(ii)] $\Tor^{R}_{i}(N,M) \in \mathcal{S}_{\mathcal{P}}(R)$ for every finitely generated $R$-module $N$ with $\Supp_{R}(N)\subseteq V(\mathfrak{a})$, and for every $i\geq 0$.
\item[(iii)] $\Tor^{R}_{i}(N,M) \in \mathcal{S}_{\mathcal{P}}(R)$ for some finitely generated $R$-module $N$ with $\Supp_{R}(N)= V(\mathfrak{a})$, and for every $0 \leq i \leq n$.
\item[(iv)] $\Ext^{i}_{R}(N,M) \in \mathcal{S}_{\mathcal{P}}(R)$ for every finitely generated $R$-module $N$ with $\Supp_{R}(N)\subseteq V(\mathfrak{a})$, and for every $i\geq 0$.
\item[(v)] $\Ext^{i}_{R}(N,M) \in \mathcal{S}_{\mathcal{P}}(R)$ for some finitely generated $R$-module $N$ with $\Supp_{R}(N)= V(\mathfrak{a})$, and for every $0 \leq i \leq n$.
\end{enumerate}
\end{corollary}

\begin{prf}
(i) $\Leftrightarrow$ (ii) and (i) $\Leftrightarrow$ (iv): Since $H_{i}(\underline{a};M)=0$ for every $i > n$, these equivalences follow from Theorems \ref{1.1}
and \ref{1.2}, respectively.

(i) $\Leftrightarrow$ (iii): Follows from Theorem \ref{1.1}.

(i) $\Leftrightarrow$ (v): Follows from Theorem \ref{1.2}.
\end{prf}

The following corollaries describe the numerical invariants $\mathcal{P}$-depth and $\mathcal{P}$-width in terms of Koszul homology, local homology, and local cohomology.

\begin{corollary} \label{3.5}
Let $\mathfrak{a}=(a_{1},\ldots,a_{n})$ be an ideal of $R$, $\underline{a}=a_{1},\ldots,a_{n}$, and $M$ an $R$-module. Let $\mathcal{P}$ be a Serre property. Then the following assertions hold:
\begin{enumerate}
\item[(i)] $\mathcal{P}$-$\depth_{R}(\mathfrak{a},M) = \inf \left\{i \geq 0 \suchthat H_{n-i}(\underline{a};M) \notin \mathcal{S}_{\mathcal{P}}(R) \right\}.$
\item[(ii)] $\mathcal{P}$-$\width_{R}(\mathfrak{a},M) = \inf \left\{i \geq 0 \suchthat H_{i}(\underline{a};M) \notin \mathcal{S}_{\mathcal{P}}(R) \right\}.$
\item[(iii)] We have $\mathcal{P}$-$\depth_{R}(\mathfrak{a},M) < \infty$ if and only if $\mathcal{P}$-$\width_{R}(\mathfrak{a},M) < \infty$. Moreover in this case, we have
\begin{center}
$\sup\left\{i \geq 0 \suchthat H_{i}(\underline{a};M) \notin \mathcal{S}_{\mathcal{P}}(R) \right\} +$ $\mathcal{P}$-$\depth_{R}(\mathfrak{a},M)=n.$
\end{center}
\end{enumerate}
\end{corollary}

\begin{prf}
(i) and (ii): Follows from Theorems \ref{1.1} and \ref{1.2}.

(iii): The first assertion follows from Corollary \ref{3.4}. For the second assertion, note that
$$\inf \left\{i \geq 0 \suchthat H_{n-i}(\underline{a};M) \notin \mathcal{S}_{\mathcal{P}}(R) \right\} = n- \sup \left\{i \geq 0 \suchthat H_{i}(\underline{a};M) \notin \mathcal{S}_{\mathcal{P}}(R) \right\}.$$
\end{prf}

\begin{corollary} \label{3.6}
Let $\mathfrak{a}$ be an ideal of $R$, and $M$ an $R$-module. Let $\mathcal{P}$ be a Serre property. Then the following assertions hold:
\begin{enumerate}
\item[(i)] If $\mathcal{P}$ satisfies the condition $\mathfrak{C}_{\mathfrak{a}}$, then
\begin{center}
$\mathcal{P}$-$\depth_{R}(\mathfrak{a},M) = \inf \left\{i \geq 0 \suchthat H^{i}_{\mathfrak{a}}(M) \notin \mathcal{S}_{\mathcal{P}}(R) \right\}.$
\end{center}
\item[(ii)] If $\mathcal{P}$ satisfies the condition $\mathfrak{D}_{\mathfrak{a}}$, then
\begin{center}
$\mathcal{P}$-$\width_{R}(\mathfrak{a},M) = \inf \left\{i \geq 0 \suchthat H^{\mathfrak{a}}_{i}(M) \notin \mathcal{S}_{\mathcal{P}}\left(\widehat{R}^{\mathfrak{a}}\right) \right\}.$
\end{center}
\item[(iii)] If $\mathcal{P}$ satisfies both conditions $\mathfrak{C}_{\mathfrak{a}}$ and $\mathfrak{D}_{\mathfrak{a}}$, and $\mathcal{P}$-$\depth_{R}(\mathfrak{a},M) < \infty$, then
\begin{center}
$\mathcal{P}$-$\depth_{R}(\mathfrak{a},M)$ $+$ $\mathcal{P}$-$\width_{R}(\mathfrak{a},M) \leq \ara(\mathfrak{a}).$
\end{center}
\end{enumerate}
\end{corollary}

\begin{prf}
(i) and (ii): Follows from Theorems \ref{1.1} and \ref{1.2}.

(iii): Clear by (i), (ii), and Corollary \ref{3.5} (iii).
\end{prf}

\section{Noetherianness and Artinianness}

In this section, we apply the results of Section 3 to obtain some characterizations of noetherian local homology modules and artinian local cohomology modules. The following result together with Corollary \ref{4.4} generalize \cite[Propositions 7.1, 7.2 and 7.4]{WW} when applied to modules.

\begin{proposition} \label{4.1}
Let $\mathfrak{a}=(a_{1},\ldots,a_{n})$ be an ideal of $R$, $\underline{a}=a_{1},\ldots,a_{n}$, and $M$ an $R$-module. Then the following assertions are equivalent for any given $s \geq 0$:
\begin{enumerate}
\item[(i)] $H_{i}(\underline{a};M)$ is a finitely generated $R$-module for every $0 \leq i \leq s$.
\item[(ii)] $\Tor^{R}_{i}(N,M)$ is a finitely generated $R$-module for every finitely generated $R$-module $N$ with $\Supp_{R}(N)\subseteq V(\mathfrak{a})$, and for every $0 \leq i \leq s$.
\item[(iii)] $\Tor^{R}_{i}(N,M)$ is a finitely generated $R$-module for some finitely generated $R$-module $N$ with $\Supp_{R}(N)= V(\mathfrak{a})$, and for every $0 \leq i \leq s$.
\item[(iv)] $H^{\mathfrak{a}}_{i}(M)$ is a finitely generated $\widehat{R}^{\mathfrak{a}}$-module for every $0 \leq i \leq s$.
\end{enumerate}
\end{proposition}

\begin{prf}
Obvious in view of Example \ref{2.6} (ii) and Theorem \ref{1.1}.
\end{prf}

The following corollary provides a characterization of noetherian local homology modules in its full generality.

\begin{corollary} \label{4.2}
Let $\mathfrak{a}=(a_{1},\ldots,a_{n})$ be an ideal of $R$, $\underline{a}=a_{1},\ldots,a_{n}$, and $M$ an $R$-module. Then the following assertions are equivalent:
\begin{enumerate}
\item[(i)] $H^{\mathfrak{a}}_{i}(M)$ is a finitely generated $\widehat{R}^{\mathfrak{a}}$-module for every $i \geq 0$.
\item[(ii)] $H_{i}(\underline{a};M)$ is a finitely generated $R$-module for every $0 \leq i \leq n$.
\item[(iii)] $H_{i}(\underline{a};M)$ is a finitely generated $R$-module for every $0 \leq i \leq \hd(\mathfrak{a},M)$.
\item[(iv)] $\Tor^{R}_{i}(N,M)$ is a finitely generated $R$-module for every finitely generated $R$-module $N$ with $\Supp_{R}(N)\subseteq V(\mathfrak{a})$, and for every $i \geq 0$.
\item[(v)] $\Tor^{R}_{i}(N,M)$ is a finitely generated $R$-module for some finitely generated $R$-module $N$ with $\Supp_{R}(N)= V(\mathfrak{a})$, and for every $0 \leq i \leq \hd(\mathfrak{a},M)$.
\item[(vi)] $\Ext^{i}_{R}(N,M)$ is a finitely generated $R$-module for every finitely generated $R$-module $N$ with $\Supp_{R}(N)\subseteq V(\mathfrak{a})$, and for every $i \geq 0$.
\item[(vii)] $\Ext^{i}_{R}(N,M)$ is a finitely generated $R$-module for some finitely generated $R$-module $N$ with $\Supp_{R}(N)= V(\mathfrak{a})$, and for every $0 \leq i \leq n$.
\end{enumerate}
\end{corollary}

\begin{prf}
(ii) $\Leftrightarrow$ (iv) $\Leftrightarrow$ (vi) $\Leftrightarrow$ (vii): Follows from Corollary \ref{3.4}.

(iii) $\Leftrightarrow$ (v): Follows from Proposition \ref{4.1} upon setting $s= \hd(\mathfrak{a},M)$.

(i) $\Leftrightarrow$ (iii): Since $H^{\mathfrak{a}}_{i}(M)=0$ for every $i > \hd(\mathfrak{a},M)$, the result follows from Proposition \ref{4.1}.

(i) $\Leftrightarrow$ (iv): Follows from Proposition \ref{4.1}.
\end{prf}

One should note that a slightly weaker version of the following result has been proved in \cite[Theorem 5.5]{Me2} by using a different method.

\begin{proposition} \label{4.3}
Let $\mathfrak{a}=(a_{1},\ldots,a_{n})$ be an ideal of $R$, $\underline{a}=a_{1},\ldots,a_{n}$, and $M$ an $R$-module. Then the following assertions are equivalent for any given $s \geq 0$:
\begin{enumerate}
\item[(i)] $H_{n-i}(\underline{a};M)$ is an artinian $R$-module for every $0 \leq i \leq s$.
\item[(ii)] $\Ext^{i}_{R}(N,M)$ is an artinian $R$-module for every finitely generated $R$-module $N$ with $\Supp_{R}(N)\subseteq V(\mathfrak{a})$, and for every $0 \leq i \leq s$.
\item[(iii)] $\Ext^{i}_{R}(N,M)$ is an artinian $R$-module for some finitely generated $R$-module $N$ with $\Supp_{R}(N)= V(\mathfrak{a})$, and for every $0 \leq i \leq s$.
\item[(iv)] $H^{i}_{\mathfrak{a}}(M)$ is an artinian $R$-module for every $0 \leq i \leq s$.
\end{enumerate}
\end{proposition}

\begin{prf}
Obvious in view of Example \ref{2.8} (ii) and Theorem \ref{1.2}.
\end{prf}

The following corollary provides a characterization of artinian local cohomology modules in its full generality.

\begin{corollary} \label{4.4}
Let $\mathfrak{a}=(a_{1},\ldots,a_{n})$ be an ideal of $R$, $\underline{a}=a_{1},\ldots,a_{n}$, and $M$ an $R$-module. Then the following assertions are equivalent:
\begin{enumerate}
\item[(i)] $H^{i}_{\mathfrak{a}}(M)$ is an artinian $R$-module for every $i \geq 0$.
\item[(ii)] $H_{i}(\underline{a};M)$ is an artinian $R$-module for every $0 \leq i \leq n$.
\item[(iii)] $H_{n-i}(\underline{a};M)$ is an artinian $R$-module for every $0 \leq i \leq \cd(\mathfrak{a},M)$.
\item[(iv)] $\Ext^{i}_{R}(N,M)$ is an artinian $R$-module for every finitely generated $R$-module $N$ with $\Supp_{R}(N)\subseteq V(\mathfrak{a})$, and for every $i \geq 0$.
\item[(v)] $\Ext^{i}_{R}(N,M)$ is an artinian $R$-module for some finitely generated $R$-module $N$ with $\Supp_{R}(N)= V(\mathfrak{a})$, and for every $0 \leq i \leq \cd(\mathfrak{a},M)$.
\item[(vi)] $\Tor^{R}_{i}(N,M)$ is an artinian $R$-module for every finitely generated $R$-module $N$ with $\Supp_{R}(N)\subseteq V(\mathfrak{a})$, and for every $i \geq 0$.
\item[(vii)] $\Tor^{R}_{i}(N,M)$ is an artinian $R$-module for some finitely generated $R$-module $N$ with $\Supp_{R}(N)= V(\mathfrak{a})$, and for every $0 \leq i \leq n$.
\end{enumerate}
\end{corollary}

\begin{prf}
(ii) $\Leftrightarrow$ (iv) $\Leftrightarrow$ (vi) $\Leftrightarrow$ (vii): Follows from Corollary \ref{3.4}.

(iii) $\Leftrightarrow$ (v): Follows from Proposition \ref{4.3} upon setting $s= \cd(\mathfrak{a},M)$.

(i) $\Leftrightarrow$ (iii): Since $H^{i}_{\mathfrak{a}}(M)=0$ for every $i > \cd(\mathfrak{a},M)$, the result follows from Proposition \ref{4.3}.

(i) $\Leftrightarrow$ (iv): Follows from Proposition \ref{4.3}.
\end{prf}

\section{Vanishing Results}

In this section, we treat the vanishing results. Since the property of being zero is a Serre property that satisfies the condition $\mathfrak{D}_{\mathfrak{a}}$, we obtain the following result.

\begin{proposition} \label{5.1}
Let $\mathfrak{a}=(a_{1},\ldots,a_{n})$ be an ideal of $R$, $\underline{a}=a_{1},\ldots,a_{n}$, and $M$ an $R$-module. Then the following assertions are equivalent for any given $s \geq 0$:
\begin{enumerate}
\item[(i)] $H_{i}(\underline{a};M)=0$ for every $0 \leq i \leq s$.
\item[(ii)] $\Tor^{R}_{i}(N,M)=0$ for every $R$-module $N$ with $\Supp_{R}(N)\subseteq V(\mathfrak{a})$, and for every $0 \leq i \leq s$.
\item[(iii)] $\Tor^{R}_{i}(N,M)=0$ for some finitely generated $R$-module $N$ with $\Supp_{R}(N)= V(\mathfrak{a})$, and for every $0 \leq i \leq s$.
\item[(iv)] $H^{\mathfrak{a}}_{i}(M)=0$ for every $0 \leq i \leq s$.
\end{enumerate}
\end{proposition}

\begin{prf}
Immediate from Theorem \ref{1.1}. For part (ii), note that every module is a direct limit of its finitely generated submodules and Tor functor commutes with direct limits.
\end{prf}

The following result is proved in \cite[Theorem 4.4]{St} using a different method, but here it is an immediate by-product of Proposition \ref{5.1}.

\begin{corollary} \label{5.2}
Let $M$ be an $R$-module, and $N$ a finitely generated $R$-module. Then the following conditions are equivalent for any given $s \geq 0$:
\begin{enumerate}
\item[(i)] $\Tor^{R}_{i}(N,M)=0$ for every $0 \leq i \leq s$.
\item[(ii)] $\Tor^{R}_{i}\left(R/ \ann_{R}(N),M\right)=0$ for every $0 \leq i \leq s$.
\end{enumerate}
\end{corollary}

\begin{prf}
Immediate from Proposition \ref{5.1}.
\end{prf}

We observe that Corollary \ref{5.3} below generalizes \cite[Corollary 4.3]{Va}, which states that if $N$ is a faithful finitely generated $R$-module, then $N\otimes_{R}M=0$ if and only if $M=0$.

\begin{corollary} \label{5.3}
Let $M$ be an $R$-module, and $N$ a finitely generated $R$-module. Then the following conditions are equivalent:
\begin{enumerate}
\item[(i)] $M\otimes_{R}N=0$.
\item[(ii)] $M= \ann_{R}(N)M$.
\end{enumerate}
In particular, we have
$$\Supp_{R}(M\otimes_{R}N)=\Supp_{R}\left(M/ \ann_{R}(N)M\right).$$
\end{corollary}

\begin{prf}
For the equivalence of (i) and (ii), let $s=0$ in Corollary \ref{5.2}. For the second part, note that given any $\mathfrak{p}\in \Spec(R)$,
we have $\ann_{R}(N)_{\frak p}=\ann_{R_{\frak p}}(N_{\frak p})$, and so $(M\otimes_{R}N)_{\mathfrak{p}}=0$ if and only if $\left(M/\ann_{R}(N)M
\right)_{\mathfrak{p}}=0.$
\end{prf}

The support formula in Corollary \ref{5.3} generalizes the well-known formula
$$\Supp_{R}(M\otimes_{R}N)=\Supp_{R}(M)\cap \Supp_{R}(N),$$
which holds whenever $M$ and $N$ are both assumed to be finitely generated.

Using Corollary \ref{3.6}, we intend to obtain two somewhat known descriptions of the numerical invariant $\width_{R}(\mathfrak{a},M)$ in terms of the Koszul homology and local homology. However, we need some generalizations of the five-term exact sequences. To the best of our knowledge, the only place where one may find such generalizations is \cite[Corollaries 10.32 and 10.34]{Ro}. But, the statements there are not correct and no proof is presented. Hence due to lack of a suitable reference, we deem it appropriate to include the correct statements with proofs for the convenience of the reader.

\begin{lemma} \label{5.4}
Let $E^{2}_{p,q} \underset {p} \Rightarrow H_{p+q}$ be a spectral sequence. Then the following assertions hold:
\begin{enumerate}
\item[(i)] If $E^{2}_{p,q} \underset {p} \Rightarrow H_{p+q}$ is first quadrant and there is an integer $n \geq 1$ such that $E^{2}_{p,q}=0$ for every $q \leq n-2$, then there is a five-term exact sequence
$$H_{n+1} \rightarrow E^{2}_{2,n-1} \rightarrow E^{2}_{0,n} \rightarrow H_{n} \rightarrow E^{2}_{1,n-1} \rightarrow 0.$$
\item[(ii)] If $E^{2}_{p,q} \underset {p} \Rightarrow H_{p+q}$ is third quadrant and there is an integer $n \geq 1$ such that $E^{2}_{p,q}=0$ for every $q \geq 2-n$, then there is a five-term exact sequence
$$0 \rightarrow E^{2}_{-1,1-n} \rightarrow H_{-n} \rightarrow E^{2}_{0,-n} \rightarrow E^{2}_{-2,1-n} \rightarrow H_{-n-1}.$$
\end{enumerate}
\end{lemma}

\begin{prf}
(i): Consider the following homomorphisms
$$0=E^{2}_{4,n-2} \xrightarrow{d^{2}_{4,n-2}} E^{2}_{2,n-1} \xrightarrow{d^{2}_{2,n-1}} E^{2}_{0,n} \xrightarrow{d^{2}_{0,n}} E^{2}_{-2,n+1}=0.$$
We thus have
$$E^{3}_{2,n-1} \cong \ker d^{2}_{2,n-1}/ \im d^{2}_{4,n-2} \cong \ker d^{2}_{2,n-1},$$
and
$$E^{3}_{0,n}\cong \ker d^{2}_{0,n}/ \im d^{2}_{2,n-1}= \coker d^{2}_{2,n-1}.$$
Let $r \geq 3$. Consider the following homomorphisms
$$0=E^{r}_{r+2,n-r} \xrightarrow{d^{r}_{r+2,n-r}} E^{r}_{2,n-1} \xrightarrow{d^{r}_{2,n-1}} E^{r}_{2-r,n+r-2}=0.$$
We thus have
$$E^{r+1}_{2,n-1} \cong \ker d^{r}_{2,n-1}/ \im d^{r}_{r+2,n-r} \cong E^{r}_{2,n-1}.$$
Therefore,
$$\ker d^{2}_{2,n-1}\cong E^{3}_{2,n-1} \cong E^{4}_{2,n-1} \cong \cdots \cong E^{\infty}_{2,n-1}.$$
Further, consider the following homomorphisms
$$0=E^{r}_{r,n-r+1} \xrightarrow{d^{r}_{r,n-r+1}} E^{r}_{0,n} \xrightarrow{d^{r}_{0,n}} E^{r}_{-r,n+r-1}=0.$$
We thus have
$$E^{r+1}_{0,n}\cong \ker d^{r}_{0,n}/ \im d^{r}_{r,n-r+1} \cong E^{r}_{0,n}.$$
Therefore,
$$\coker d^{2}_{2,n-1} \cong E^{3}_{0,n} \cong E^{4}_{0,n} \cong \cdots \cong E^{\infty}_{0,n}.$$
Hence, we get the following exact sequence
\begin{equation} \label{eq:5.4.1}
0 \rightarrow E^{\infty}_{2,n-1} \rightarrow E^{2}_{2,n-1} \xrightarrow{d^{2}_{2,n-1}} E^{2}_{0,n} \rightarrow E^{\infty}_{0,n} \rightarrow 0.
\end{equation}
There is a finite filtration
$$0=U^{-1} \subseteq U^{0} \subseteq \cdots \subseteq U^{n+1}=H_{n+1},$$
such that $E^{\infty}_{p,n+1-p}\cong U^{p}/U^{p-1}$ for every $0 \leq p \leq n+1$. If $3 \leq p \leq n+1$, then $n+1-p \leq n-2$, so $E^{\infty}_{p,n+1-p}=0$. It follows that
$$U^{2}=U^{3}=\cdots =U^{n+1}.$$
Now since $E^{\infty}_{2,n-1}\cong U^{2}/U^{1}$, we get an exact sequence
$$H_{n+1}\rightarrow E^{\infty}_{2,n-1}\rightarrow 0.$$
Splicing this exact sequence to the exact sequence \eqref{eq:5.4.1}, we get the following exact sequence
\begin{equation} \label{eq:5.4.2}
H_{n+1} \rightarrow E^{2}_{2,n-1} \xrightarrow{d^{2}_{2,n-1}} E^{2}_{0,n} \rightarrow E^{\infty}_{0,n} \rightarrow 0.
\end{equation}
On the other hand, there is a finite filtration
$$0=V^{-1} \subseteq V^{0} \subseteq \cdots \subseteq V^{n}=H_{n},$$
such that $E^{\infty}_{p,n-p}\cong V^{p}/V^{p-1}$ for every $0 \leq p \leq n$. If $2 \leq p \leq n$, then $n-p \leq n-2$, so $E^{\infty}_{p,n-p}=0$. It follows that
$$V^{1}=V^{2}=\cdots =V^{n}.$$
As $E^{\infty}_{0,n}\cong V^{0}/V^{-1} = V^{0}$ and $E^{\infty}_{1,n-1} \cong V^{1}/V^{0}$, we get the short exact sequence
$$0 \rightarrow E^{\infty}_{0,n} \rightarrow H_{n} \rightarrow E^{\infty}_{1,n-1} \rightarrow 0.$$
Splicing this short exact sequence to the exact sequence \eqref{eq:5.4.2}, yields the exact sequence
\begin{equation} \label{eq:5.4.3}
H_{n+1} \rightarrow E^{2}_{2,n-1} \xrightarrow{d^{2}_{2,n-1}} E^{2}_{0,n} \rightarrow H_{n} \rightarrow E^{\infty}_{1,n-1} \rightarrow 0.
\end{equation}
Let $r \geq 2$, and consider the following homomorphisms
$$0=E^{r}_{r+1,n-r} \xrightarrow{d^{r}_{r+1,n-r}} E^{r}_{1,n-1} \xrightarrow{d^{r}_{1,n-1}} E^{r}_{1-r,n+r-2}=0.$$
We thus have
$$E^{r+1}_{1,n-1}\cong \ker d^{r}_{1,n-1}/ \im d^{r}_{r+1,n-r} \cong E^{r}_{1,n-1}.$$
Therefore,
$$E^{2}_{1,n-1}\cong E^{3}_{1,n-1}\cong \cdots \cong E^{\infty}_{1,n-1}.$$
Thus from the exact sequence \eqref{eq:5.4.3}, we get the desired exact sequence
$$H_{n+1} \rightarrow E^{2}_{2,n-1} \xrightarrow{d^{2}_{2,n-1}} E^{2}_{0,n} \rightarrow H_{n} \rightarrow E^{2}_{1,n-1} \rightarrow 0.$$

(ii): Similar to (i).
\end{prf}

For the next result, we need to recall the notion of coassociated prime ideals. Given an $R$-module $M$, a prime ideal $\mathfrak{p}\in \Spec (R)$ is said to be a coassociated prime ideal of $M$ if $\mathfrak{p}= \ann_{R}(M/N)$ for some submodule $N$ of $M$ such that $M/N$ is an artinian $R$-module. The set of coassociated prime ideals of $M$ is denoted by $\Coass_{R}(M)$.

\begin{corollary} \label{5.5}
Let $\mathfrak{a}=(a_{1},\ldots,a_{n})$ be an ideal of $R$, $\underline{a}=a_{1},\ldots,a_{n}$, and $M$ an $R$-module. Then the following assertions hold:
\begin{enumerate}
\item[(i)] $\width_{R}(\mathfrak{a},M) = \inf \left\{i \geq 0 \suchthat H_{i}(\underline{a};M) \neq 0 \right\} = \inf \left\{i \geq 0 \suchthat H_{i}^{\mathfrak{a}}(M) \neq 0 \right\}$.
\item[(ii)] $\Tor^{R}_{\width_{R}(\mathfrak{a},M)}(R/\mathfrak{a},M) \cong (R/\mathfrak{a})\otimes_{R}H_{\width_{R}(\mathfrak{a},M)}^{\mathfrak{a}}(M)$.
\item[(iii)] $\Lambda^{\mathfrak{a}}\left(H_{\width_{R}(\mathfrak{a},M)}^{\mathfrak{a}}(M)\right) \cong \underset{n}\varprojlim \Tor^{R}_{\width_{R}(\mathfrak{a},M)}(R/\mathfrak{a}^{n},M)$.
\item[(iv)] $\Coass_{R}\left(\Tor^{R}_{\width_{R}(\mathfrak{a},M)}(R/\mathfrak{a},M)\right)=\Coass_{R}\left(H_{\width_{R}(\mathfrak{a},M)}^{\mathfrak{a}}(M)\right) \cap V(\mathfrak{a})$.
\end{enumerate}
\end{corollary}

\begin{prf}
(i): Follows immediately from Corollaries \ref{3.5} and \ref{3.6}.

(ii): Consider the first quadrant spectral sequence
$$E^{2}_{p,q}= \Tor^{R}_{p}\left(R/\mathfrak{a},H^{\mathfrak{a}}_{q}(M)\right) \underset {p} \Rightarrow \Tor^{R}_{p+q}(R/\mathfrak{a},M).$$
Let $n=\width_{R}(\mathfrak{a},M)$. Then by (i), $E^{2}_{p,q}=0$ for every $q\leq n-1$. Now, Lemma \ref{5.4} (i) gives the exact sequence
$$\Tor^{R}_{n+1}(R/\mathfrak{a},M) \rightarrow \Tor^{R}_{2}\left(R/\mathfrak{a},H^{\mathfrak{a}}_{n-1}(M)\right) \rightarrow (R/\mathfrak{a})
\otimes_{R} H^{\mathfrak{a}}_{n}(M)\rightarrow $$$$ \Tor^{R}_{n}(R/\mathfrak{a},M) \rightarrow \Tor_{1}^{R}\left(R/\mathfrak{a},
H^{\mathfrak{a}}_{n-1}(M)\right) \rightarrow 0.$$
But
$$\Tor^{R}_{2}\left(R/\mathfrak{a},H^{\mathfrak{a}}_{n-1}(M)\right)=0=\Tor_{1}^{R}\left(R/\mathfrak{a},H^{\mathfrak{a}}_{n-1}(M)\right),$$
so
$$(R/\mathfrak{a})\otimes_{R} H^{\mathfrak{a}}_{n}(M) \cong \Tor^{R}_{n}(R/\mathfrak{a},M).$$

(iii): Using (ii) and the facts that $H_{i}^{\mathfrak{a}^{n}}(M)\cong H_{i}^{\mathfrak{a}}(M)$ and $\width_{R}(\mathfrak{a}^{n},M)=\width_{R}(\mathfrak{a},M)$
for every $n,i\geq 0$, we get:
\begin{equation*}
\begin{split}
\Lambda^{\mathfrak{a}}\left(H_{\width_{R}(\mathfrak{a},M)}^{\mathfrak{a}}(M)\right) & = \underset{n}\varprojlim \left(H_{\width_{R}(\mathfrak{a},M)}^{\mathfrak{a}}(M)/\mathfrak{a}^{n}H_{\width_{R}(\mathfrak{a},M)}^{\mathfrak{a}}(M)\right) \\
 & \cong \underset{n}\varprojlim \Tor^{R}_{\width_{R}(\mathfrak{a},M)}(R/\mathfrak{a}^{n},M).
\end{split}
\end{equation*}

(iv): Using (ii) and \cite[Theorem 1.21]{Ya}, we have
\begin{equation*}
\begin{split}
\Coass_{R}\left(\Tor^{R}_{\width_{R}(\mathfrak{a},M)}(R/\mathfrak{a},M)\right) & = \Coass_{R}\left((R/\mathfrak{a})\otimes_{R}H_{\width_{R}(\mathfrak{a},M)}^{\mathfrak{a}}(M)\right) \\
 & = \Supp_{R}(R/\mathfrak{a}) \cap \Coass_{R}\left(H_{\width_{R}(\mathfrak{a},M)}^{\mathfrak{a}}(M)\right) \\
 & = V(\mathfrak{a}) \cap \Coass_{R}\left(H_{\width_{R}(\mathfrak{a},M)}^{\mathfrak{a}}(M)\right).
\end{split}
\end{equation*}
\end{prf}

Note that part (iii) of Corollary \ref{5.5} is proved in \cite[Proposition 2.5]{Si2} by deploying a different method. On the other hand, in parallel with Corollary \ref{5.10} (iii) below, one may wonder if intersecting with $V(\mathfrak{a})$ in part (iv) of Corollary \ref{5.5} is redundant. In other words, $\Coass_{R}\left(H_{\width_{R}(\mathfrak{a},M)}^{\mathfrak{a}}(M)\right)$ may be contained in $V(\mathfrak{a})$. However, the following example shows that this is not the case in general.

\begin{example} \label{5.6}
Let $R:=\mathbb{Q}[X,Y]_{(X,Y)}$ and $\mathfrak{m}:=(X,Y)R$. Then $\width_{R}(\mathfrak{m},R)=0$ and
$$H^{\frak m}_0(R)\cong \widehat{R}^{\mathfrak{m}}\cong \mathbb{Q}[[X,Y]].$$
For each $n \in \mathbb{Z}$, let $\fp_n:=(X-nY)R$. Then it is easy to see that $R/\fp_n\cong \mathbb{Q}[Y]_{(Y)}$, and so it is not a complete local ring. By \cite [Beispiel 2.4]{Z},
\begin{equation*}
\begin{split}
\Coass_R\left(H^{\frak m}_0(R)\right) & = \Coass_R\left(\widehat{R}^{\mathfrak{m}}\right) \\
 & = \{\fm\} \cup \left\{\mathfrak{p}\in \Spec R \suchthat R/\mathfrak{p} \text{ is not a complete local ring} \right\}.
\end{split}
\end{equation*}
Hence $\Coass_{R}\left(H_0^{\frak m}(R)\right)$ is not a finite set, while $\Coass_{R}\left(H_0^{\frak m}(R)\right)\cap V(\mathfrak{m})= \{\mathfrak{m}\}$. In particular, $\Coass_{R}\left(H_0^{\frak m}(R)\right) \nsubseteq V(\mathfrak{m})$.
\end{example}

Since the property of being zero is a Serre property that satisfies the condition $\mathfrak{C}_{\mathfrak{a}}$, we obtain the following result.

\begin{proposition} \label{5.7}
Let $\mathfrak{a}=(a_{1},\ldots,a_{n})$ be an ideal of $R$, $\underline{a}=a_{1},\ldots,a_{n}$, and $M$ an $R$-module. Then the following assertions are equivalent for any given $s \geq 0$:
\begin{enumerate}
\item[(i)] $H_{n-i}(\underline{a};M)=0$ for every $0 \leq i \leq s$.
\item[(ii)] $\Ext^{i}_{R}(N,M)=0$ for every finitely generated $R$-module $N$ with $\Supp_{R}(N)\subseteq V(\mathfrak{a})$, and for every $0 \leq i \leq s$.
\item[(iii)] $\Ext^{i}_{R}(N,M)=0$ for some finitely generated $R$-module $N$ with $\Supp_{R}(N)= V(\mathfrak{a})$, and for every $0 \leq i \leq s$.
\item[(iv)] $H^{i}_{\mathfrak{a}}(M)=0$ for every $0 \leq i \leq s$.
\end{enumerate}
\end{proposition}

\begin{prf}
Immediate from Theorem \ref{1.2}.
\end{prf}

The following result is proved in \cite[Theorem 3.2]{St} using a different method, but here it is an immediate by-product of Proposition \ref{5.7}.

\begin{corollary} \label{5.8}
Let $M$ be an $R$-module, and $N$ a finitely generated $R$-module. Then the following conditions are equivalent for any given $s \geq 0$:
\begin{enumerate}
\item[(i)] $\Ext^{i}_{R}(N,M)=0$ for every $0 \leq i \leq s$.
\item[(ii)] $\Ext^{i}_{R}\left(R/ \ann_{R}(N),M\right)=0$ for every $0 \leq i \leq s$.
\end{enumerate}
\end{corollary}

\begin{prf}
Immediate from Proposition \ref{5.7}.
\end{prf}

The following special case may be of independent interest.

\begin{corollary} \label{5.9}
Let $M$ be an $R$-module, and $N$ a finitely generated $R$-module. Then the following conditions are equivalent:
\begin{enumerate}
\item[(i)] $\Hom_{R}(N,M)=0$.
\item[(ii)] $\left(0:_{M} \ann_{R}(N)\right)=0$.
\end{enumerate}
\end{corollary}

\begin{prf}
Let $s=0$ in Corollary \ref{5.8}.
\end{prf}

It is easy to deduce from Corollary \ref{5.9} that given a finitely generated $R$-module $N$, we have $\Hom_{R}(N,M)\neq 0$ if and only if there are elements $x\in N$ and $0\neq y \in M$ with $\ann_{R}(x) \subseteq \ann_{R}(y)$, which is known as the Hom Vanishing Lemma in \cite[Page 11]{C}.

We state the dual result to Corollary \ref{5.5} for the sake of integrity and completeness. Parts (ii) and (iii) of Corollary \ref{5.10} below are stated in \cite[Proposition 1.1]{Mar}, and it is only mentioned that part (ii) can be deduced from a spectral sequence. In addition, a proof of this result is offered in \cite[Corollary 2.3]{Me1} by using different techniques.

\begin{corollary} \label{5.10}
Let $\mathfrak{a}=(a_{1},\ldots,a_{n})$ be an ideal of $R$, $\underline{a}=a_{1},\ldots,a_{n}$, and $M$ an $R$-module. Then the following assertions hold:
\begin{enumerate}
\item[(i)] $\depth_{R}(\mathfrak{a},M) = \inf \left\{i \geq 0 \suchthat H_{n-i}(\underline{a};M) \neq 0 \right\} = \inf \left\{i \geq 0 \suchthat H_{\mathfrak{a}}^{i}(M) \neq 0 \right\}$.
\item[(ii)] $\Ext^{\depth_{R}(\mathfrak{a},M)}_{R}(R/\mathfrak{a},M) \cong \Hom_{R}\left(R/\mathfrak{a},H^{\depth_{R}(\mathfrak{a},M)}_{\mathfrak{a}}(M)\right)$.
\item[(iii)] $\Ass_{R}\left(\Ext^{\depth_{R}(\mathfrak{a},M)}_{R}(R/\mathfrak{a},M)\right)=\Ass_{R}\left(H^{\depth_{R}(\mathfrak{a},M)}_{\mathfrak{a}}(M)\right)$.
\end{enumerate}
\end{corollary}

\begin{prf}
(i): Follows immediately from Corollaries \ref{3.5} and \ref{3.6}.

(ii): Consider the third quadrant spectral sequence
$$E^{2}_{p,q}= \Ext^{-p}_{R}\left(R/\mathfrak{a},H^{-q}_{\mathfrak{a}}(M)\right) \underset {p} \Rightarrow \Ext^{-p-q}_{R}(R/\mathfrak{a},M),$$
and use Lemma \ref{5.4} (ii).

(iii): Follows from (ii).
\end{prf}

Finally, we present the following comprehensive vanishing result.

\begin{corollary} \label{5.11}
Let $\mathfrak{a}=(a_{1},\ldots,a_{n})$ be an ideal of $R$, $\underline{a}=a_{1},\ldots,a_{n}$, and $M$ an $R$-module. Then the following assertions are equivalent:
\begin{enumerate}
\item[(i)] $H_{i}(\underline{a};M)=0$ for every $0 \leq i \leq n$.
\item[(ii)] $H_{i}(\underline{a};M)=0$ for every $0 \leq i \leq \hd(\mathfrak{a},M)$.
\item[(iii)] $H_{i}(\underline{a};M)=0$ for every $n-\cd(\mathfrak{a},M) \leq i \leq n$.
\item[(iv)] $\Tor^{R}_{i}(N,M)=0$ for every $R$-module $N$ with $\Supp_{R}(N)\subseteq V(\mathfrak{a})$, and for every $i \geq 0$.
\item[(v)] $\Tor^{R}_{i}(N,M)=0$ for some finitely generated $R$-module $N$ with $\Supp_{R}(N)= V(\mathfrak{a})$, and for every $0 \leq i \leq \hd(\mathfrak{a},M)$.
\item[(vi)] $\Ext^{i}_{R}(N,M)=0$ for every finitely generated $R$-module $N$ with $\Supp_{R}(N)\subseteq V(\mathfrak{a})$, and for every $i \geq 0$.
\item[(vii)] $\Ext^{i}_{R}(N,M)=0$ for some finitely generated $R$-module $N$ with $\Supp_{R}(N)= V(\mathfrak{a})$, and for every $0 \leq i \leq \cd(\mathfrak{a},M)$.
\item[(viii)] $H^{\mathfrak{a}}_{i}(M)=0$ for every $i \geq 0$.
\item[(ix)] $H^{i}_{\mathfrak{a}}(M)=0$ for every $i \geq 0$.
\end{enumerate}
\end{corollary}

\begin{prf}
Follows from Propositions \ref{5.1} and \ref{5.7}.
\end{prf}

The following corollary is proved in \cite[Corollary 1.7]{Si2}. However, it is an immediate consequence of the results thus far obtained.

\begin{corollary} \label{5.12}
Let $\mathfrak{a}$ be an ideal of $R$, and $M$ an $R$-module. Then $\depth_{R}(\mathfrak{a},M) < \infty$ if and only if $\width_{R}(\mathfrak{a},M) < \infty$.
Moreover in this case, we have $$\depth_{R}(\mathfrak{a},M)+ \width_{R}(\mathfrak{a},M) \leq \ara(\mathfrak{a}).$$
\end{corollary}

\begin{prf}
Clear by Corollary \ref{3.6} (iii).
\end{prf}


\end{document}